\RequirePackage{ifpdf}
\ifpdf 
\documentclass[pdftex]{sigma}
\else
\documentclass{sigma}
\fi

\begin{document}

\allowdisplaybreaks

\renewcommand{\thefootnote}{$\star$}

\renewcommand{\PaperNumber}{010}

\FirstPageHeading

\ShortArticleName{The Rational qKZ Equation and Shifted Non-Symmetric Jack Polynomials}

\ArticleName{The Rational qKZ Equation\\ and Shifted Non-Symmetric Jack Polynomials\footnote{This paper is a contribution to the Special
Issue on Dunkl Operators and Related Topics. The full collection
is available at
\href{http://www.emis.de/journals/SIGMA/Dunkl_operators.html}{http://www.emis.de/journals/SIGMA/Dunkl\_{}operators.html}}}


\Author{Saburo KAKEI~$^\dag$, Michitomo NISHIZAWA~$^\ddag$,
 Yoshihisa SAITO~$^\mathsection$ and Yoshihiro TAKEYAMA~$^\mathparagraph$}

\AuthorNameForHeading{S. Kakei, M. Nishizawa, Y. Saito and Y. Takeyama}

\Address{$^\dag$~Department of Mathematics, College of Science, Rikkyo University,\\
\hphantom{$^\dag$}~Nishi-Ikebukuro, Toshima-ku, Tokyo 171-8501, Japan}
\EmailD{\href{mailto:kakei@rikkyo.ac.jp}{kakei@rikkyo.ac.jp}}

\Address{$^\ddag$~Department of Mathematics,
Faculty of Education, Hirosaki University,\\
\hphantom{$^\ddag$}~1 Bunkyo-cho, Hirosaki, Aomori 036-8560, Japan}
\EmailD{\href{mailto:mnishi@cc.hirosaki-u.ac.jp}{mnishi@cc.hirosaki-u.ac.jp}}

\Address{$^\mathsection$~Graduate School of Mathematical Sciences, University of Tokyo,
Tokyo 153-8914, Japan}
\EmailD{\href{mailto:yosihisa@ms.u-tokyo.ac.jp}{yosihisa@ms.u-tokyo.ac.jp}}

\Address{$^\mathparagraph$~Department of Mathematics,
Graduate School of Pure and Applied Sciences,\\
\hphantom{$^\mathparagraph$}~University of Tsukuba, Tsukuba, Ibaraki 305-8571, Japan}
\EmailD{\href{mailto:takeyama@math.tsukuba.ac.jp}{takeyama@math.tsukuba.ac.jp}}

\ArticleDates{Received October 15, 2008, in f\/inal form January 15,
2009; Published online January 27, 2009}

\Abstract{We construct special solutions to the rational quantum Knizhnik--Zamolodchikov equation
associated with the Lie algebra $gl_{N}$.
The main ingredient is a special class of the shifted non-symmetric Jack polynomials.
It may be regarded as a shifted version of the singular polynomials studied by Dunkl.
We prove that our solutions contain those obtained as
a scaling limit of matrix elements of the vertex operators of level one.}


\Classification{39A13; 33C52; 81R50}

\renewcommand{\thefootnote}{\arabic{footnote}}
\setcounter{footnote}{0}

\section{Introduction}

The quantum Knizhnik--Zamolodchikov (qKZ) equation,
derived by Frenkel and Reshetikhin~\cite{FR}, is a system of
dif\/ference equations satisf\/ied by matrix elements of vertex operators
in the representation theory of the quantum af\/f\/ine algebras.
In this paper we consider the rational version of the qKZ equation.
Let $V=\oplus_{\epsilon=1}^{N}\mathbb{C}v_{\epsilon}$ be the vector representation
of the Lie algebra~$gl_{N}$.
The rational qKZ equation is
the following system of dif\/ference equations for an unknown function
$G(z_{1}, \ldots , z_{n})$ taking values in $V^{\otimes n}$:
\begin{gather*}
 G(z_{1}, \ldots , z_{m}-\kappa, \ldots , z_{n})=
R_{m,m-1}(z_{m}-z_{m-1}-\kappa) \cdots R_{m,1}(z_{m}-z_{1}-\kappa) \\
 \qquad {}\times
\Bigg(  \prod_{j=1}^{N}p_{j}^{h_{j}} \Bigg)_{m}
R_{m,n}(z_{m}-z_{n}) \cdots R_{m,m+1}(z_{m}-z_{m+1})
G(z_{1}, \ldots , z_{m}, \ldots , z_{n})
\end{gather*}
for $m=1, \ldots , n$.
Here $h_{i}$ $(1 \le i \le N)$ is a basis of the Cartan subalgebra,
$p_{1}, \ldots , p_{N}$ and $\kappa$ are nonzero parameters,
and $R(z) \in {\rm End}(V^{\otimes 2})$ is the rational $R$-matrix (see \eqref{eq:R-matrix} below).
The lower indices of the operators signify the components of $V^{\otimes n}$ on which they act.
The value $-(N+\kappa)$ is called the level of the equation.

In \cite{KT} special solutions to the qKZ equation are constructed.
Let us recall the construction brief\/ly.
Consider the coef\/f\/icients in an expansion of the unknown function
with respect to the standard basis of $V^{\otimes n}$.
In \cite{KT} a suf\/f\/icient condition is given so that
they form a solution to the qKZ equation.
It is described as an eigenvalue problem for the $q$-Dunkl--Cherednik operators
and the Demazure--Lusztig operators.
The problem can be solved by using the non-symmetric Macdonald polynomials~\cite{C1, M},
which are by def\/inition eigenfunctions of the $q$-Dunkl--Cherednik operators.
Once one obtains a solution to the eigenvalue problem,
the coef\/f\/icients are generated from it
through the action of the Demazure--Lusztig operators.
Thus we obtain special solutions to the qKZ equation from
the non-symmetric Macdonald polynomials.
In this construction the polynomial representation of the af\/f\/ine Hecke algebra
is a main tool.

In this article we apply the method of~\cite{KT} to the rational qKZ equation.
The main tool is the degenerate double af\/f\/ine Hecke algebra $\widetilde{\bf H}_{\kappa}$.
It is a trigonometric degeneration of the double af\/f\/ine Hecke algebra
introduced by Cherednik~\cite{C}.
Usually we consider the polynomial representation of $\widetilde{\bf H}_{\kappa}$
given by the dif\/ferential operators called the Dunkl--Cherednik operators.
They form a commuting family of dif\/ferential operators and
the polynomial eigenfunctions are the non-symmetric Jack polynomials \cite{O}.
In our construction of solutions to the rational qKZ equation,
we consider another polynomial representation given in terms of
commuting dif\/ference operators $Y_{i}$ $(1 \le i \le n)$, which are dif\/ference analogues of
the Dunkl operators
\begin{gather*}
\mathcal{D}_{i}:=
\kappa \frac{\partial}{\partial x_{i}}+\sum_{j(\not=i)}\frac{1}{x_{i}-x_{j}}(1-\sigma_{ij}),
\end{gather*}
where $\sigma_{ij}$ is the transposition of variables $x_{i} \leftrightarrow x_{j}$.
The parameter $\kappa$ will play the role of the dif\/ference step in the rational qKZ equation.
The operator $Y_{i}$ has an expansion $Y_{i}=1-\mathcal{D}_{i}+\widetilde{Y}_{i}$,
where $\widetilde{Y}_{i}$ decreases the total degree by at least two.
We can construct another family of commuting dif\/ference operators
in the same way as the construction of the Dunkl--Cherednik operators from the Dunkl operators.
Then the corresponding polynomial eigenfunctions are
the shifted non-symmetric Jack polynomials,
which are originally obtained by Knop~\cite{Knop} as a~classical limit
of the non-symmetric quantum Capelli polynomials.

In our construction of special solutions to the rational qKZ equation
we need a polynomial eigenfunction of the operators $Y_{i}$ $(1 \le i \le n)$.
Then, from the expansion $Y_{i}=1-\mathcal{D}_{i}+\widetilde{Y}_{i}$,
the eigenvalue should be one,
and the leading degree part of a polynomial eigenfunction
should belong to $\cap_{i=1}^{n}\ker{\mathcal{D}_{i}}$
(note that any eigenvalue of the Dunkl operator is equal to zero).
A~nonzero polynomial in the joint kernel is called a singular polynomial,
which is deeply studied by Dunkl~\mbox{\cite{D1, D2}}.
It is proved that constants are only singular polynomials for generic $\kappa$,
and specif\/ic values of $\kappa$, called singular values,
such that there exists a non-constant singular polynomial
are completely determined by Dunkl, de Jeu and Opdam
in more general setting associated with f\/inite ref\/lection groups~\cite{DDO}.
When $\kappa$ is equal to a singular value,
the space of singular polynomials becomes a non-trivial $\mathbb{C}\mathfrak{S}_{n}$-module.
In~\cite{D1, D2} Dunkl determined its irreducible components and
proved that all the multiplicities are equal to one.
Each irreducible component has a basis consisting of non-symmetric Jack polynomials
which are well-def\/ined at the singular value and singular.
In this paper we consider the shifted case.
We call a polynomial eigenfunction of $Y_{i}$ $(1 \le i \le n)$
a {\it shifted singular polynomial}.
In Theorem \ref{thm:singular-iso} below we prove that
the space of shifted singular polynomials is isomorphic to
that of singular polynomials as $\mathbb{C}\mathfrak{S}_{n}$-modules,
and it has a basis consisting of shifted non-symmetric Jack polynomials.
This fact follows from the inversion formula \cite{Knop}, which relates
the non-symmetric Jack polynomials and the shifted ones.
We can get a solution to the rational qKZ equation from
a shifted singular polynomial which satisf\/ies some symmetry conditions
(see \eqref{eq:eigenvalue-problem-2} below).

Now we mention the result \cite{DZ} by Di Francesco and Zinn-Justin.
They constructed a class of polynomial solutions to the qKZ equation and
observed that their scaling limit give ``extended Joseph polynomials''.
As mentioned in~\cite{DZ}, the scaling limit of their solutions solves
the rational qKZ equation of level one.
It had been used in the construction of a new integral formula for solutions to
the rational qKZ equation of level zero in~\cite{T}.
We prove that our solutions contain the special solutions above of level one.
In our framework, the extended Joseph polynomials appearing in
the solutions to the rational qKZ equation are shifted non-symmetric Jack polynomials
which are singular.
The leading degree part of the extended Joseph polynomials is nothing but
the Joseph polynomial~\cite{J} which is a singular non-symmetric Jack polynomial.
Thus the inversion formula for (shifted) non-symmetric Jack polynomials at a singular value
may be regarded as a map which bridges some data of cohomology
of orbital varieties (Joseph polynomials)
and their equivariant counterpart (extended Joseph polynomials).

The rest of the paper is organized as follows.
In Section~\ref{sec:shifted-Jack} we summarize the def\/initions and properties of
(shifted) non-symmetric Jack polynomials from a viewpoint of
the representation theory of the degenerate double af\/f\/ine Hecke algebra.
In Section~\ref{sec:rational-qKZ} we give the rational qKZ equation
and a corresponding eigenvalue problem.
In Section~\ref{sec:solution} we construct polynomial solutions to the eigenvalue problem,
and hence we obtain special solutions to the rational qKZ equation.
The key is the isomorphism as $\mathbb{C}\mathfrak{S}_{n}$-modules
between the space of singular polynomials and that of shifted singular polynomials.
At the very last we prove an explicit formula for
some singular (shifted) non-symmetric Jack polynomials.

\section{Shifted non-symmetric Jack polynomials}\label{sec:shifted-Jack}

We introduce some notation.
We call an element of $\mathbb{Z}_{\ge 0}^{n}$ {\it a composition}.
In the following we f\/ix $n$, and  denote by $\Lambda$ the set of compositions.
A composition $\lambda=(\lambda_{1}, \ldots , \lambda_{n})$ is called {\it dominant}
if $\lambda_{1} \ge \cdots \ge \lambda_{n}$.
The symmetric group $\mathfrak{S}_{n}$ acts on $\Lambda$ by
$\sigma \lambda:=(\lambda_{\sigma^{-1}(1)}, \ldots , \lambda_{\sigma^{-1}(n)})$.
For $\lambda \in \Lambda$, denote by $\lambda^{+}$
the unique dominant element in the orbit $\mathfrak{S}_{n}\lambda$,
and by $w_{\lambda}^{+}$ the shortest element of $\mathfrak{S}_{n}$ such that
$w_{\lambda}^{+}\lambda^{+}{=}\lambda$.
We def\/ine the dominance order $\le$ on $\Lambda$ by
$\lambda \le \mu \!\Leftrightarrow\!  \sum\limits_{i=1}^{k}\lambda_{i}\le \sum\limits_{i=1}^{k}\mu_{i}$~$(\forall\,{k})$.

\subsection{Non-symmetric Jack polynomials}

First let us recall the def\/inition of the degenerate double af\/f\/ine Hecke algebra
(or the trigonometric Cherednik algebra) following \cite{Suzuki}\footnote{We change the notation in \cite{Suzuki} as
$x_{i} \to x_{n+1-i}$, $s_{i} \to s_{n-i}$, $u_{i} \to u_{n+1-i}$.}.
Let $\kappa$ be an indeterminate and
denote the coef\/f\/icient f\/ield by $\mathbb{F}=\mathbb{C}(\kappa)$.

\begin{definition}\label{def:dDAHA}
The {\it degenerate double affine Hecke algebra} $\widetilde{\bf H}_{\kappa}$
is the unital associative $\mathbb{F}$-algebra
generated by
\begin{gather*}
x_{i}^{\pm 1} \quad (1 \le i \le n), \qquad
s_{i} \quad (1 \le i <n), \qquad
u_{i} \quad (1 \le i \le n)
\end{gather*}
satisfying the following relations: 
\begin{gather*}
[x_{i}, x_{j}]=0, \qquad [u_{i}, u_{j}]=0, \\
s_{i}^2=1, \qquad s_{i}s_{i+1}s_{i}=s_{i+1}s_{i}s_{i+1}, \qquad
s_{i}s_{j}=s_{j}s_{i} \qquad (|i-j|>1), \\
s_{i}x_{i}s_{i}=x_{i+1}, \qquad s_{i}x_{j}=x_{j}s_{i} \qquad (j \not=i, i+1), \\
s_{i}u_{i}=u_{i+1}s_{i}+1, \qquad s_{i}u_{j}=u_{j}s_{i} \qquad (j\not=i, i+1), \\
[u_{i}, x_{j}]=\left\{
\begin{array}{ll}
{}-x_{i}s_{ji} & (j<i), \\
\kappa x_{j}+x_{j}\sum\limits_{1\le k<j}s_{kj}+\sum\limits_{j<k\le n}s_{jk}x_{j} & (j=i),  \\
{}-x_{j}s_{ij} & (i<j).
\end{array}
\right.
\end{gather*}
Here $s_{ij}=(s_{i} \cdots s_{j-1})(s_{j-2} \cdots s_{i})$ for $i<j$. 
\end{definition}

Note that $\widetilde{\bf H}_{\kappa}$ contains
the group algebra of the symmetric group
$\langle s_{1}, \ldots , s_{n-1} \rangle \simeq \mathbb{F}\mathfrak{S}_{n}$.

Denote by $\widetilde{\bf H}_{\kappa}^{+}$ the subalgebra
generated by $x_{i}$, $u_{i}$ $(1 \le i \le n)$ and $s_{i}$ $(1 \le i <n)$.
The algebra~$\widetilde{\bf H}_{\kappa}^{+}$
and its polynomial representation are crucial tools
in the theory of Jack polynomials.
Let $\mathbb{F}[x]=\mathbb{F}[x_{1}, \ldots , x_{n}]$ be the polynomial ring.
Denote by $\sigma_{i,j}$ the transposition of variables
$x_{i} \leftrightarrow x_{j}$.
Then the polynomial representation
$\phi : \widetilde{\bf H}_{\kappa}^{+} \to {\rm End}\,\mathbb{F}[x]$ is def\/ined by
\begin{gather}
\phi(x_{i})=x_{i}, \qquad
\phi(s_{i})=\sigma_{i,i+1}, \qquad
\phi(u_{i})=x_{i}\mathcal{D}_{i}+\sum_{i<k\le n}\sigma_{ik},
\label{eq:Dunkl-Cherednik}
\end{gather}
where $\mathcal{D}_{i}$ is the {\it Dunkl operator}
\begin{gather}
\mathcal{D}_{i}:=
\kappa \frac{\partial}{\partial x_{i}}+\sum_{j(\not=i)}\frac{1}{x_{i}-x_{j}}(1-\sigma_{ij}).
\label{eq:Dunkl}
\end{gather}
The operators $\phi(u_{i})$ $(1\le i \le n)$ are called
the {\it Dunkl--Cherednik operators}.

The {\it non-symmetric Jack polynomials} $\xi_{\lambda}$ $(\lambda \in \Lambda)$
are uniquely def\/ined by the following conditions:
\begin{gather*}
\phi(u_{i})\xi_{\lambda}=c_{i}(\lambda)\xi_{\lambda} \qquad (1\le i \le n), \qquad
\xi_{\lambda}(x)=x^{\lambda}+\sum_{\mu < \lambda}c_{\lambda\mu}x^{\mu} \qquad
(c_{\lambda\mu} \in \mathbb{F}).
\end{gather*}
Here $x^{\lambda}=x_{1}^{\lambda_{1}} \cdots x_{n}^{\lambda_{n}}$, and
$c(\lambda)=(c_{1}(\lambda), \ldots , c_{n}(\lambda)) \in \mathbb{F}^{n}$ is def\/ined by
\begin{gather*}
c(\lambda):=\kappa \lambda+w_{\lambda}^{+}\rho,
\end{gather*}
where $\rho:=(n-1, \ldots , 1, 0)$.

\subsection{Shifted non-symmetric Jack polynomials}

The shifted non-symmetric Jack polynomials are def\/ined by
the vanishing conditions.
Hereafter we use the variables $z=(z_{1}, \ldots , z_{n})$ for
the shifted non-symmetric Jack polynomials,
and denote by $\mathbb{F}[z]=\mathbb{F}[z_{1}, \ldots , z_{n}]$
the ring of polynomials in $z_{1}, \ldots , z_{n}$.

\begin{theorem}[\protect{\cite[Theorem 6.2]{Knop}}]\label{thm:def-shifted-nonsymmetric-Jack}
For $\lambda \in \Lambda$ there exists
a polynomial $E_{\lambda}(z) \in \mathbb{F}[z]$
uniquely determined by the following conditions:
\begin{gather*}
E_{\lambda}(c(\mu))=0 \qquad \hbox{for any $\mu \in \Lambda$ such that} \quad
|\mu|\le |\lambda|, \  \mu\not=\lambda, \\
E_{\lambda}(z)=\sum_{|\nu|\le |\lambda|}d_{\lambda\nu}z^{\nu} \qquad
(d_{\lambda\nu} \in \mathbb{F}, \  d_{\lambda\lambda}=1).
\end{gather*}
Here $|\lambda|:=\sum\limits_{i=1}^{n}\lambda_{i}$.
We call $E_{\lambda}$ $(\lambda \in \Lambda)$
the shifted non-symmetric Jack polynomials.
\end{theorem}

\begin{remark}\label{rem:Knop}
The shifted non-symmetric Jack polynomials def\/ined in Theorem~\ref{thm:def-shifted-nonsymmetric-Jack}
are slightly dif\/ferent from the original ones given by Knop~\cite{Knop} as follows.
In \cite{Knop} the parameter is denoted by $r$.
Let $\widetilde{\mathbb{F}}:=\mathbb{C}(r)$ be the corresponding coef\/f\/icient f\/ield.
Def\/ine a $\mathbb{C}$-algebra isomorphism
$\Psi: \mathbb{F}[z] \to \widetilde{\mathbb{F}}[z]$ by
$\Psi(\kappa)=1/r$ and $\Psi(z_{i})=n-1+z_{i}/r$.
Then the polynomials $\widetilde{E}_{\lambda}$ def\/ined in~\cite[Theorem 6.2]{Knop} are
related to $E_{\lambda}$ in Theorem~\ref{thm:def-shifted-nonsymmetric-Jack}
by $\Psi(E_{\lambda})=r^{-|\lambda|}\widetilde{E}_{\lambda}$.
\end{remark}

The shifted non-symmetric Jack polynomials are also characterized
as joint eigenfunctions of some commuting dif\/ference operators.
To construct the operators,
we introduce another polynomial representation of $\widetilde{\bf H}_{\kappa}$.
Def\/ine the operators $T_{i}$ $(1 \le i <n)$ and $\omega$ acting on $\mathbb{F}[z]$ by
\begin{gather}
T_{i}:=\sigma_{i,i+1}+\frac{1}{z_{i}-z_{i+1}}(1-\sigma_{i,i+1}),
\label{eq:Demazure-Lusztig} \\
(\omega f)(z_{1}, \ldots , z_{n}):=f(z_{n}-\kappa, z_{1}, \ldots , z_{n-1}).
\label{eq:Dynkin-auto}
\end{gather}
Here $\sigma_{i,i+1}$ is the transposition of variables $z_{i} \leftrightarrow z_{i+1}$.
Set
\begin{gather}
Y_{i}:=T_{i} \cdots T_{n-1}   \omega
T_{1} \cdots T_{i-1} \qquad
(1 \le i \le n).
\label{eq:q-Dunkl}
\end{gather}

\begin{proposition}\label{prop:polynomial-representation}
There exists an $\mathbb{F}$-algebra homomorphism
$\psi : \widetilde{\bf H}_{\kappa} \to {\rm End}\,\mathbb{F}[z]$ such that
\begin{gather*}
\psi(x_{i})=Y_{i}, \qquad
\psi(s_{i})=T_{i}, \qquad
\psi(u_{i})=z_{i}.
\end{gather*}
\end{proposition}

\begin{proof}
It can be checked by direct calculation using
$\omega T_{i}=T_{i-1}\omega$ $(1<i<n)$.
\end{proof}

Note that the generator $x_{i}$ acts as a dif\/ference operator
and $u_{i}$ as the multiplication of the variable $z_{i}$ on $\mathbb{F}[z]$.
It is opposite to the action $\phi$ \eqref{eq:Dunkl-Cherednik}.

\begin{remark}
The usual polynomial representation $\phi$ is induced from
the one-dimensional representation of the subalgebra generated by
$s_{i}$ $(1 \le i <n)$ and $u_{i}$ $(1 \le i \le n)$.
The representation~$\psi$ in Proposition~\ref{prop:polynomial-representation}
is constructed from an opposite side in the following sense.
Consider the subalgebra generated by $s_{i}$ $(1 \le i <n)$ and $x_{i}$ $(1 \le i \le n)$.
It is isomorphic to the group algebra of the extended af\/f\/ine Weyl group
$\widehat{\mathfrak{S}}_{n}$.
Then the induced module of the trivial $\mathbb{C}\widehat{\mathfrak{S}}_{n}$-module
is isomorphic to $\mathbb{F}[u_{1}, \ldots , u_{n}]$ as an $\mathbb{F}$-vector space.
Change the variables $u_{i} \to z_{i}$ and we obtain the polynomial representation $\psi$.
\end{remark}

\begin{proposition}\label{prop:commuting-v}
Define $v_{i} \in \widetilde{\bf H}_{\kappa}$ $(1 \le i \le n)$ by
\begin{gather*}
v_{i}:=u_{i}(1-x_{i})+\sum_{i<k \le n}s_{ik}x_{i}.
\end{gather*}
Then the following relations hold:
\begin{gather*}
[v_{i}, v_{j}]=0, \qquad
s_{i}v_{i}=v_{i+1}s_{i}+1, \qquad
s_{i}v_{j}=v_{j}s_{i} \qquad (j\not=i, i+1).
\end{gather*}
\end{proposition}

\begin{proof}
Set $\widehat{u}_{i}:=u_{i}-v_{i}=\big(u_{i}-\sum\limits_{i<k\le n}s_{ik}\big)x_{i}$.
Then we can check that
$[\widehat{u}_{i}, \widehat{u}_{j}]=0$ and
$[v_{j}, \widehat{u}_{i}]=[v_{i}, \widehat{u}_{j}]=-\widehat{u}_{j}s_{ij}$ for $i<j$.
Thus we get the f\/irst relation.
The others can be proved by easy calculation.
\end{proof}

Thus we obtain commuting dif\/ference operators $\psi(v_{i})$ $(1 \le i \le n)$ acting on $\mathbb{F}[z]$.
The shifted non-symmetric Jack polynomials are joint eigenfunctions of them.

\begin{theorem}[\protect{\cite[Theorem 6.6]{Knop}}]
The shifted non-symmetric Jack polynomials $E_{\lambda}$ $(\lambda \in \Lambda)$
are uniquely characterized by the following properties:
\begin{gather*}
\psi(v_{i})E_{\lambda}=c_{i}(\lambda)E_{\lambda} \qquad (1 \le i \le n), \\
E_{\lambda}(z)=\sum_{|\nu|\le |\lambda|}d_{\lambda\nu}z^{\nu} \qquad
(d_{\lambda\nu} \in \mathbb{F}, \  d_{\lambda\lambda}=1).
\end{gather*}
\end{theorem}

\begin{remark}
The operators $\psi(v_{i})$ $(1 \le i \le n)$ are related to $\widetilde{\Xi}_{i}$'s
in~\cite[Theorem~6.6]{Knop}   by $\Psi^{-1}\widetilde{\Xi}_{i}\Psi=(\psi(v_{i})-(n-1))/\kappa$,
where $\Psi$ is the map def\/ined in Remark~\ref{rem:Knop}.
\end{remark}

Here we brief\/ly explain why the eigenvalues of $\psi(v_{i})$ $(1 \le i \le n)$
are the same as those of the Dunkl--Cherednik operators $\phi(u_{i})$ $(1 \le i \le n)$.
Since the operator $T_{i}-\sigma_{i}$ decreases the total degree by one,
we have an expansion $Y_{i}=1-\mathcal{D}_{i}+\widetilde{Y}_{i}$,
where $\mathcal{D}_{i}$ is the Dunkl operator~\eqref{eq:Dunkl},
and $\widetilde{Y}_{i}$ is an operator which decreases the total degree by at least two.
Thus we see that $\psi(v_{i})-\phi(u_{i})$ is an operator
which strictly decreases the total degree.
Hence $\psi(v_{i})$ is triangular with respect to the order on monomials in $\mathbb{F}[z]$
induced from the dominance order $\le$ and the total degree,
and its eigenvalues are equal to those of $\phi(u_{i})$.

\begin{proposition}\label{prop:eigenvalue-Y}
Any eigenvalue of the operators $Y_{i}$ $(1 \le i \le n)$ on~$\mathbb{F}[z]$
is equal to one.
\end{proposition}

\begin{proof}
{}From the argument above,
the operator $Y_{i}-1=-\mathcal{D}_{i}+\widetilde{Y}_{i}$ decreases the total degree by at least one.
This implies the proposition.
\end{proof}

\subsection{The inversion formula}


The $\widetilde{\bf H}_{\kappa}$-module $\mathbb{F}[z]$ is naturally
an $\widetilde{\bf H}_{\kappa}^{+}$-module by restriction.
However we consider another $\widetilde{\bf H}_{\kappa}^{+}$-module structure
def\/ined through the following homomorphism.

\begin{proposition}[\cite{F}]
Set $\widehat{u}_{i}:=u_{i}-v_{i}=\big(u_{i}-\sum\limits_{i<k\le n}s_{ik}\big)x_{i}$ $(1 \le i \le n)$.
Then there exists an $\mathbb{F}$-algebra homomorphism $\tau$ satisfying
\begin{gather*}
\tau(x_{i})=\widehat{u}_{i}, \qquad \tau(s_{i})=s_{i}, \qquad \tau(u_{i})=v_{i}.
\end{gather*}
\end{proposition}

Hereafter we regard $\mathbb{F}[z]$ as an $\widetilde{\bf H}_{\kappa}^{+}$-module
through the map $\psi \circ \tau$.
To clarify it we write ${}^{\tau}\mathbb{F}[z]$ instead of $\mathbb{F}[z]$
when we consider the $\widetilde{\bf H}_{\kappa}^{+}$-module structure.
On the other hand we regard $\mathbb{F}[x]$ as an $\widetilde{\bf H}_{\kappa}^{+}$-module
by the restriction $\phi|_{\widetilde{\bf H}_{\kappa}^{+}}$.

\begin{proposition}
The map
\begin{gather}
\Phi : \ \ \mathbb{F}[x] \longrightarrow {}^{\tau}\mathbb{F}[z], \qquad
f(x_{1}, \ldots , x_{n}) \mapsto
\psi(f(\widehat{u}_{1}, \ldots , \widehat{u}_{n}))\,1
\label{eq:inversion-generic}
\end{gather}
is an isomorphism between $\widetilde{\bf H}_{\kappa}^{+}$-modules.
The inverse $\Phi^{-1}$ is a map taking the leading degree part and changing the variables $z_{i} \to x_{i}$.
\end{proposition}

\begin{proof}
The intertwining property follows from
$\phi(s_{i})\,1=1$, $\phi(u_{i})\,1=n-i$ in $\mathbb{F}[x]$, and
$\psi(\tau(s_{i}))\,1=\psi(s_{i})\,1=1$ and $\psi(\tau(u_{i}))\,1=\psi(v_{i})\,1=n-i$ in
${}^{\tau}\mathbb{F}[z]$.
The operator $\psi(\widehat{u}_{i})-z_{i}$ decreases the total degree
of polynomials in $z$, and hence
$\Phi(x^{\lambda})=z^{\lambda}+(\hbox{lower degree terms})$ for any $\lambda \in \Lambda$.
Therefore $\Phi$ is an isomorphism.
\end{proof}

\begin{corollary}[\protect{\cite[Theorem 6.9]{Knop}}]\label{cor:inversion}
$\Phi(\xi_{\lambda})=E_{\lambda}$ for any $\lambda \in \Lambda$.
Hence the leading degree part of $E_{\lambda}$ is equal to $\xi_{\lambda}$.
\end{corollary}

\begin{proof}
Apply $\Phi$ to the both hand sides of $\phi(u_{i}) \xi_{\lambda}=c_{i}(\lambda) \xi_{\lambda}$.
Then we see that $\Phi(\xi_{\lambda})$ is an eigenfunction of
$\psi(v_{i})$ with the eigenvalue $c_{i}(\lambda)$.
The coef\/f\/icient of $z^{\lambda}$ in $\Phi(\xi_{\lambda})$ is equal to one.
This implies that $\Phi(\xi_{\lambda})=E_{\lambda}$.
\end{proof}

\begin{remark}
The operators $\psi(\widehat{u}_{i})$ $(1 \le i \le n)$ are related to $\widetilde{Z}_{i}$'s
in~\cite[Theorem~6.9]{Knop}  by $\Psi^{-1}\widetilde{Z}_{i}\Psi=\psi(\widehat{u}_{i})/\kappa$,
where $\Psi$ is the map def\/ined in Remark~\ref{rem:Knop}.
\end{remark}

The equality $\Phi(\xi_{\lambda})=E_{\lambda}$ is called the {\it inversion formula}.

The non-symmetric Jack polynomials $\xi_{\lambda}$ may have poles
when $\kappa$ is a non-positive rational number.
From the inversion formula, we see that the shifted one $E_{\lambda}$ also has the same property.
We say that $\xi_{\lambda}$ (or $E_{\lambda}$) is
{\it well-def\/ined at $\kappa=\kappa_{0} \in \mathbb{Q}_{\le 0}$}
if it does not have a pole at $\kappa=\kappa_{0}$.

\begin{corollary}\label{cor:well-defined}
For $\lambda \in \Lambda$, $\xi_{\lambda}$ is well-defined at $\kappa=\kappa_{0} \in \mathbb{Q}_{\le 0}$
if and only if $E_{\lambda}$ is well-defined.
\end{corollary}

In Section \ref{sec:special-solutions}
we specialize $\kappa$ to some constant $\kappa=\kappa_{0} \in \mathbb{Q}^{\times}$.
Then the $\mathbb{C}$-algebra $\widetilde{\bf H}_{\kappa_{0}}^{+}$ acts
on $\mathbb{C}[x]$ and ${}^{\tau}\mathbb{C}[z]$
through the maps $\phi$ and $\psi\circ \tau$, respectively.
In this situation the map
\begin{gather}
\Phi: \ \ \mathbb{C}[x] \longrightarrow {}^{\tau}\mathbb{C}[z], \qquad
f(x_{1}, \ldots , x_{n}) \mapsto \psi(f(\widehat{u}_{1}, \ldots , \widehat{u}_{n}))\,1
\label{eq:def-Phi}
\end{gather}
is still an isomorphism between $\widetilde{\bf H}_{\kappa_{0}}^{+}$-modules.

\section{The rational qKZ equation}\label{sec:rational-qKZ}

\subsection{The rational qKZ equation}

Let $V=\oplus_{\epsilon=1}^{N}\mathbb{C}v_{\epsilon}$ be the $N$-dimensional vector space.
We regard $V$ as the vector representation of the Lie algebra $gl_{N}$.
Denote by $h_{i}$ $(1 \le i \le N)$ a basis of the Cartan subalgebra,
which acts on $V$ by $h_{i}v_{\epsilon}=\delta_{i,\epsilon}v_{\epsilon}$.

The rational $R$-matrix $R(z) \in {\rm End}(V^{\otimes 2})$ is def\/ined by
\begin{gather}
R(z)=\frac{z+P}{z+1},
\label{eq:R-matrix}
\end{gather}
where $P$ is the transposition: $P(u \otimes v)=v\otimes u$.

Fix a parameter $\kappa \in \mathbb{C}^{\times}$.
The {\it rational quantum Knizhnik--Zamolodchikov (rational qKZ) equation} is
the following system of dif\/ference equations for an unknown function
$G(z_{1}, \ldots , z_{n})$ taking values in $V^{\otimes n}$:
\begin{gather*}
G(z_{1}, \ldots , z_{m}-\kappa, \ldots , z_{n})=
R_{m,m-1}(z_{m}-z_{m-1}-\kappa) \cdots R_{m,1}(z_{m}-z_{1}-\kappa) \\
\qquad{}\times
\Bigg(   \prod_{j=1}^{N}p_{j}^{h_{j}} \Bigg)_{m}
R_{m,n}(z_{m}-z_{n}) \cdots R_{m,m+1}(z_{m}-z_{m+1})
G(z_{1}, \ldots , z_{m}, \ldots , z_{n})
\end{gather*}
for $m=1, \ldots , n$.
Here $p_{1}, \ldots , p_{N}$ are nonzero parameters,
and the lower indices of the operators signify the components of $V^{\otimes n}$ on which they act.
The value $-(N+\kappa)$ is called the {\it level} of the equation.

\subsection{Eigenvalue problem associated with the rational qKZ equation}

Hereafter we assume that
\begin{gather*}
n \ge N \ge 2.
\end{gather*}

Fix positive integers $d_{1}, \ldots , d_{N}$ satisfying $\sum\limits_{j=1}^{N}d_{j}=n$,
and def\/ine an $n$-tuple $\delta$ by
$\delta=(1^{d_{1}}, \ldots , N^{d_{N}})$.
Consider the following eigenvalue problem
for a function $F_{\delta}=F_{\delta}(z_{1}, \ldots , z_{n})$:
\begin{gather}
Y_{i}F_{\delta}  =\chi_{i}F_{\delta} \qquad
(1 \le i \le n, \, \chi_{i}\in \mathbb{C}^{\times}),
\label{eq:eigenvalue-problem-1}
\\
T_{i}F_{\delta}  =-F_{\delta} \qquad {\rm if} \quad \delta_{i}=\delta_{i+1}.
\label{eq:eigenvalue-problem-2}
\end{gather}
Here $Y_{i}$ and $T_{i}$ are given by~\eqref{eq:q-Dunkl} and~\eqref{eq:Demazure-Lusztig}, respectively.
Note that the consistency of~\eqref{eq:eigenvalue-problem-1} and~\eqref{eq:eigenvalue-problem-2}
implies that $\chi_{i}=\chi_{i+1}$ if $\delta_{i}=\delta_{i+1}$.
Hence there are $N$ independent parameters
$\chi_{d_{1}+\cdots +d_{j}}$  $(1\le j \le N)$.

Once one solves the eigenvalue problem above,
a solution to the rational qKZ equation can be constructed as follows.
Suppose that $F_{\delta}$ is a solution to the eigenvalue problem~\eqref{eq:eigenvalue-problem-1} and~\eqref{eq:eigenvalue-problem-2}.
Set $I_{d_{1}, \ldots , d_{N}}:=\mathfrak{S}_{n}\delta \subset \{1, \ldots , N\}^{n}$.
For $\epsilon=(\epsilon_{1}, \ldots , \epsilon_{n}) \in I_{d_{1}, \ldots , d_{N}}$,
take a shortest element $w_{\epsilon} \in \mathfrak{S}_{n}$ such that
$w_{\epsilon}\epsilon=\delta$.
Let $w_{\epsilon}=s_{i_{1}} \cdots s_{i_{l}}$ be a reduced expression,
and set $T_{w_{\epsilon}}:=T_{i_{1}} \cdots T_{i_{l}}$.
This does not depend on a choice of the reduced expression.
Now def\/ine a function $F_{\epsilon}(z_{1}, \ldots , z_{n})$ by
\begin{gather*}
F_{\epsilon}:=(-1)^{\ell(w_{\epsilon})}T_{w_{\epsilon}}F_{\delta},
\end{gather*}
where $\ell(w)$ is the length of $w \in \mathfrak{S}_{n}$.
Set
\begin{gather*}
F(z_{1}, \ldots , z_{n}):=\sum_{\epsilon \in I_{d_{1}, \ldots , d_{N}}}
F_{\epsilon}(z_{1}, \ldots , z_{n})\,
v_{\epsilon_{1}} \otimes \cdots \otimes v_{\epsilon_{n}}.
\end{gather*}
Let $K(z_{1}, \ldots, z_{n})$ be a function satisfying
\begin{gather}
\frac{K(\ldots , z_{m}-\kappa, \ldots)}{K(\ldots , z_{m}, \ldots)}=
\prod_{j=1}^{m-1}\frac{z_{j}-z_{m}+\kappa+1}{z_{j}-z_{m}+\kappa-1}
\prod_{j=m+1}^{n}\frac{z_{m}-z_{j}-1}{z_{m}-z_{j}+1}
\label{eq:relation-for-K}
\end{gather}
for $m=1, \ldots , n$.
For example, take
\begin{gather}
K(z_{1}, \ldots, z_{n})=\prod_{1\le i<j \le n}
\frac{\Gamma\left(-\frac{z_{i}-z_{j}-1}{\kappa}\right)}{\Gamma\left(-\frac{z_{i}-z_{j}+1}{\kappa}\right)}.
\label{eq:K-Gamma}
\end{gather}

\begin{theorem}\label{thm:qKZ}
Set $G(z_{1}, \ldots , z_{n}):=K(z_{1}, \ldots , z_{n})F(z_{1}, \ldots , z_{n})$.
Then $G$ is a solution of the rational qKZ equation with
$p_{j}=\chi_{d_{1}+\cdots +d_{j}}$ $(1\le j \le N)$.
\end{theorem}

\begin{proof}
In the same way as the proof of Theorem 3.6 in~\cite{KT},
we can see that the following relations hold:
\begin{gather*}
T_{i}F_{\ldots , \epsilon_{i}, \epsilon_{i+1}, \ldots }
 =
{}-F_{\ldots , \epsilon_{i+1}, \epsilon_{i}, \ldots },
\\
\omega F_{\epsilon_{n}, \epsilon_{1}, \ldots , \epsilon_{n-1}}
 =
(-1)^{n-1}\chi_{d_{1}+\cdots +d_{\epsilon_{n}}}
F_{\epsilon_{1}, \ldots , \epsilon_{n}},
\end{gather*}
where $\omega$ is def\/ined by \eqref{eq:Dynkin-auto}.
By setting $p_{j}=\chi_{d_{1}+\cdots +d_{j}}$, the two relations above
are equivalent to
\begin{gather*}
\frac{z_{i}-z_{i+1}+1}{z_{i+1}-z_{i}+1}P_{i,i+1}R_{i,i+1}(z_{i}-z_{i+1})
F(\ldots , z_{i}, z_{i+1}, \ldots)=F(\ldots , z_{i+1}, z_{i}, \ldots), \\
P_{n-1,n}\cdots P_{1,2}F(z_{n}-\kappa, z_{1}, \ldots , z_{n-1})=(-1)^{n-1}
\left(\prod_{j=1}^{N}p_{j}^{h_{j}}\right)_{n}   F(z_{1}, \ldots , z_{n}).
\end{gather*}
Combining them and \eqref{eq:relation-for-K},
we see that $G=KF$ is a solution to the rational qKZ equation.
\end{proof}

\section{Special solutions to the rational qKZ equation}\label{sec:solution}

Now we construct a class of polynomial solutions to the eigenvalue problem
\eqref{eq:eigenvalue-problem-1} and \eqref{eq:eigenvalue-problem-2}.
They create solutions to the rational qKZ equation as discussed in the previous subsection
(see Theo\-rem~\ref{thm:qKZ}).

\subsection{Shifted singular polynomials}\label{sec:shifted-singular}

First let us consider the equation \eqref{eq:eigenvalue-problem-1}.
We want to obtain polynomial solutions.
Then the eigenvalues $\chi_{i}$ $(1 \le i \le n)$ should be one
{}from Proposition \ref{prop:eigenvalue-Y}.
As seen below the polynomial eigenfunctions can be regarded as
a shifted version of singular polynomials.

A polynomial $f(x_{1}, \ldots , x_{n}) \in \mathbb{C}[x]$ is called a {\it singular polynomial}
if $f \in \cap_{i=1}^{n}\ker{\mathcal{D}_{i}}$ \cite{D1, D2, DDO}.
Denote by $\mathcal{S}$ the subspace of $\mathbb{C}[x]$ consisting of singular polynomials.
Recall that $\mathbb{C}[x]$ has $\widetilde{\bf H}_{\kappa}$-module structure
def\/ined by the map $\phi$ \eqref{eq:Dunkl-Cherednik}.
{}From the commutation relations
$\sigma_{i} \mathcal{D}_{i}=\mathcal{D}_{i+1}\sigma_{i}$ and
$\sigma_{i} \mathcal{D}_{j}=\mathcal{D}_{j} \sigma_{i}$ $(j\not=i, i+1)$,
the maps $\phi(s_{i})$ $(1 \le i <n)$ preserve $\mathcal{S}$.
Thus $\mathcal{S}$ is a $\mathbb{C}\mathfrak{S}_{n}$-module.

Now let us consider the shifted version.

\begin{definition}
We call $g(z_{1}, \ldots , z_{n}) \in \mathbb{C}[z]$ a {\it shifted singular polynomial}
if it satisf\/ies $Y_{i}g=g$ for all $1 \le i \le n$.
\end{definition}

Denote by $\widetilde{\mathcal{S}}$ the subspace of $\mathbb{C}[z]$ consisting
of shifted singular polynomials.
{}From the relations $T_{i}Y_{i}=Y_{i+1}T_{i}$
and $T_{i}Y_{j}=Y_{j}T_{i}$ $(j\not=i, i+1)$,
which follow from Proposition~\ref{prop:polynomial-representation},
the maps~$T_{i}$ $(1 \le i <n)$ preserve $\widetilde{\mathcal{S}}$.
Hence $\widetilde{\mathcal{S}}$ is also a $\mathbb{C}\mathfrak{S}_{n}$-module.
It is naturally obtained from the $\widetilde{\bf H}_{\kappa}^{+}$-module structure on
${}^{\tau}\mathbb{C}[z]$ (note that $\tau(s_{i})=s_{i}$ for all $i$).

\begin{theorem}\label{thm:singular-iso}
The map $\Phi$ \eqref{eq:def-Phi} gives a $\mathbb{C}\mathfrak{S}_{n}$-module isomorphism
$\Phi|_{\mathcal{S}}: \mathcal{S} \simeq \widetilde{\mathcal{S}}$.
\end{theorem}

\begin{proof}
Recall that $\Phi^{-1}$ is the map taking the leading degree part and
changing the variables $z_{i} \to x_{i}$.
Then we easily see $\Phi^{-1}(\widetilde{\mathcal{S}}) \subset \mathcal{S}$
 from the expansion $Y_{i}=1-\mathcal{D}_{i}+\widetilde{Y}_{i}$
given in the argument preceding to Proposition~\ref{prop:eigenvalue-Y}.

Let us prove $\Phi(\mathcal{S}) \subset \widetilde{\mathcal{S}}$.
Suppose that $f \in \mathbb{C}[x]$ be a singular polynomial.
Then we have $x_{i}\mathcal{D}_{i}f=0$ $(1 \le i \le n)$.
{}From $x_{i}\mathcal{D}_{i}=\phi(u_{i}-\sum\limits_{i<k\le n}s_{ik})$ and
the intertwining property of $\Phi$, we have
\begin{gather*}
0 =\Phi(x_{i}\mathcal{D}_{i}f)=\psi\Bigg(\tau\Bigg(u_{i}-\sum_{i<k\le n}s_{ik}\Bigg)\Bigg)\Phi(f) \\
\phantom{0}
{}=\psi(v_{i}-\sum_{i<k\le n}s_{ik})\Phi(f)=
\Bigg(z_{i}-\sum_{i<k\le n}T_{ik}\Bigg)(1-Y_{i})\Phi(f),
\end{gather*}
where $T_{ik}=\psi(s_{ik})$.
Since the operator $z_{i}-\sum\limits_{i<k\le n}T_{ik}$ is injective,
we get $Y_{i}\Phi(f)=\Phi(f)$ for $1 \le i \le n$.
\end{proof}

If $\kappa$ is generic, we have $\mathcal{S}=\mathbb{C}$, that is,
constants are only singular polynomials.
A specif\/ic parameter value $\kappa=\kappa_{0}$ is called a {\it singular value}
if there exists a non-constant singular polynomial at $\kappa=\kappa_{0}$.
In \cite{DDO} it is proved that the set of singular values is equal to
$\{-l/m \, | \, m \in \mathbb{Z}_{>0}$, $l=2, \ldots , n \,\, \hbox{and} \,\, m \not\in l\mathbb{Z} \}$\footnote{Note that the parameter $\kappa$ in \cite{D1, D2} is equal to $1/\kappa$ in this article.
See the def\/inition of the Dunkl operator~\eqref{eq:Dunkl}.}.
From Corollary~\ref{cor:inversion}, Corollary~\ref{cor:well-defined}
and Theorem~\ref{thm:singular-iso}, we f\/ind

\begin{corollary}\label{cor:singular-basis-correspondence}
Let $\kappa_{0}$ be a singular value.
Suppose that $\lambda$ is a composition such that
$\xi_{\lambda}$ is well-defined and singular at $\kappa=\kappa_{0}$.
Then $E_{\lambda}$ is also well-defined and shifted singular.
\end{corollary}

In \cite{D1, D2} Dunkl obtained an explicit description of $\mathcal{S}$
as a $\mathbb{C}\mathfrak{S}_{n}$-module.
At a singular value $\kappa=\kappa_{0}$,
the module $\mathcal{S}$ is multiplicity-free and
the set of the irreducible components is completely determined from
$\kappa_{0}$ and $n$.
Each irreducible component has a basis consisting of
non-symmetric Jack polynomials which are well-def\/ined and singular.
From Theorem~\ref{thm:singular-iso} and Corollary~\ref{cor:singular-basis-correspondence},
the $\mathbb{C}\mathfrak{S}_{n}$-module structure of $\widetilde{\mathcal{S}}$
is the same as~$\mathcal{S}$.
In the shifted case, $\widetilde{\mathcal{S}}$ has a basis
consisting of shifted non-symmetric Jack polynomials.

Now we give a limiting procedure from non-symmetric Macdonald polynomials to
shifted singular ones.
Let us recall the def\/inition of non-symmetric Macdonald polynomials.
Let~$q$ and~$t^{1/2}$ be indeterminates and
denote the coef\/f\/icient f\/ield by $\mathbb{K}=\mathbb{C}(q, t^{1/2})$.
Consider the following operators $\widehat{T}_{i}$ $(1\le i<n)$ and $\widehat{\omega}$
acting on the polynomial ring $\mathbb{K}[X_{1}, \ldots , X_{n}]$:
\begin{gather*}
\widehat{T}_{i}:=t^{1/2}\sigma_{i}+\frac{t^{1/2}-t^{-1/2}}{X_{i}/X_{i+1}-1}(\sigma_{i}-1), \qquad
(\widehat{\omega}f)(X_{1}, \ldots , X_{n}):=f(qX_{n}, X_{1}, \ldots , X_{n-1}),
\end{gather*}
where $\sigma_{i}$ is the transposition of variables $X_{i} \leftrightarrow X_{i+1}$.
The operators $\widehat{T}_{i}$ $(1 \le i<n)$ are called the {\it Demazure--Lusztig operators}.
Then the {\it $q$-Dunkl--Cherednik operators} $\widehat{Y}_{i}$ $(1\le i\le n)$ are
def\/ined by
\begin{gather*}
\widehat{Y}_{i}:=\widehat{T}_{i} \cdots \widehat{T}_{n-1} \, \widehat{\omega}\,
\widehat{T}_{1}^{-1} \cdots \widehat{T}_{i-1}^{-1}.
\end{gather*}
The {\it non-symmetric Macdonald polynomials} $E_{\lambda}^{\,q,t}$ $(\lambda \in \Lambda)$ are
def\/ined by the following properties:{\samepage
\begin{gather*}
 \widehat{Y}_{i}E_{\lambda}^{q, t}=t^{(w_{\lambda}^{+}\rho)_{i}}q^{\lambda_{i}}\,E_{\lambda}^{\,q, t}
\qquad (1 \le i \le n), \qquad
E_{\lambda}^{\,q, t}(X)=X^{\lambda}+\sum_{\mu <\lambda}a_{\mu\lambda}X^{\mu} \qquad
(a_{\mu\lambda}\in \mathbb{K}).
\end{gather*}
Here $X^{\lambda}:=X_{1}^{\lambda_{1}} \cdots X_{n}^{\lambda_{n}}$.}

Introduce a small parameter $\epsilon$ and set
\begin{gather}
X_{j}=e^{\epsilon z_{j}}, \qquad
t^{1/2}=e^{-\epsilon/2}, \qquad
q=e^{-\kappa \epsilon}.
\label{eq:limit}
\end{gather}
Then we f\/ind $\widehat{T}_{i}^{\pm 1} \to T_{i}$ and $\widehat{\omega} \to \omega$
in the limit $\epsilon \to 0$,
and hence $\widehat{Y}_{i} \to Y_{i}$.
Suppose that  $E_{\lambda}^{q, t}$
has the  Laurent expansion at $\epsilon=0$ in the form
\begin{gather*}
E_{\lambda}^{q, t}(X)|_{q=t^{\kappa}}=\epsilon^{\alpha}(g_{\lambda}(z)+o(1))
\end{gather*}
for some $\alpha \in \mathbb{Z}$ and nonzero $g_{\lambda}(z) \in \mathbb{F}[z]$.
Then we obtain $Y_{i}g_{\lambda}=g_{\lambda}$ for $1 \le i \le n$,
that is, $g_{\lambda}$ is shifted singular.
{}From the discussion preceding to Corollary~\ref{cor:singular-basis-correspondence},
$g_{\lambda}$ can become non-constant only if $\kappa$ is a singular value.
In Proposition~\ref{prop:factorized-formula} below, we give an example
of $\lambda$ and a singular value $\kappa_{0}$
such that $g_{\lambda}$ is non-constant at $\kappa=\kappa_{0}$.

\subsection{Construction of special solutions}\label{sec:special-solutions}

To construct a solution to the rational qKZ equation,
we need a shifted singular polynomial satisfying
the condition \eqref{eq:eigenvalue-problem-2}.

Let $r$ and $k$ be positive integers such that $r \ge 2$, $k \le n$ and
${\rm gcd}(r-1, k+1)=1$.
Hereafter we set\footnote{The parametrization \eqref{eq:singular-value} is motivated by
 the condition of {\it $(k,r)$-admissibility} employed in the construction
 of special solutions to the (trigonometric) qKZ equation~\cite{KT}.
 It is originally appeared in~\cite{FJMM}.}
\begin{gather}
\kappa=-\frac{k+1}{r-1}.
\label{eq:singular-value}
\end{gather}

Take three integers $s, m$ and $l$ satisfying
\begin{gather*}
n=(s+m)k+s+l, \qquad
s \ge 0, \qquad m \ge 1, \qquad 0 \le l<k.
\end{gather*}
For $a \ge 1$, we set
\begin{gather*}
\lambda[r,s,a]:=((r-1)(s+a), (r-1)(s+a-1), \ldots , (r-1)(s+1), 0).
\end{gather*}
For $a=0$, set $\lambda[r,s,0]=(0)$.
Now we def\/ine the composition $\lambda(k,s,m,l,r)\in \Lambda$ by
\begin{gather*}
\lambda(k,s,m,l,r):=\big(0^{(k+1)s}, \lambda[r,s,m-1]^{k-l}, \lambda[r,s,m]^{l}\big).
\end{gather*}

Next we take an $N$-tuple $(d_{1}, \ldots , d_{N})$ satisfying the following conditions.
For $1 \le i \le (k+1)s$, we have $d_{i}=1$.
The rest $(d_{(k+1)s+1}, \ldots , d_{N})$ should be a subdivision of $(m^{k-l}, (m+1)^{l})$,
that is, $d_{i}>0$ and
\begin{gather*}
d_{i_{j}}+\cdots +d_{i_{j+1}-1}=\left\{
\begin{array}{ll}
m & (1 \le j \le k-l), \\
m+1 & (k-l<j\le k)
\end{array}
\right.
\end{gather*}
for some $i_{1}=(k+1)s+1<i_{2}<\cdots <i_{k}<i_{k+1}=N+1$.

Now set $\delta=(1^{d_{1}}, \ldots , N^{d_{N}})$.
Then we f\/ind

\begin{proposition}\label{prop:polynomial-solution}
$E_{\lambda(k,s,m,l,r)}$ is well-defined at \eqref{eq:singular-value}
and satisfies \eqref{eq:eigenvalue-problem-1} and~\eqref{eq:eigenvalue-problem-2}
with $\chi_{i}=1$ $(1 \le i \le n)$.
Therefore the shifted non-symmetric Jack polynomial $E_{\lambda(k,s,m,l,r)}$
at \eqref{eq:singular-value} creates a~solution to the rational qKZ equation
with $p_{j}=1$ $(1 \le j \le N)$ of level $\frac{k+1}{r-1}-N$.
\end{proposition}

\begin{proof}
In \cite[Theorem 5.7]{D1} it is proved that $\xi_{\lambda(k,s,m,l,r)}$ is well-def\/ined
and singular at \eqref{eq:singular-value}
(note that if $n=k$ then $(s,m,l)=(0,1,0)$ and $E_{\lambda(k,0,1,0,r)}=E_{0^{k}}=1$ is
trivially well-def\/ined and singular).
Hence, from Corollary~\ref{cor:singular-basis-correspondence},
$E_{\lambda(k,s,m,l,r)}$ is also well-def\/ined and satisf\/ies~\eqref{eq:eigenvalue-problem-1}.

To prove that $E_{\lambda(k,s,m,l,r)}$ satisf\/ies \eqref{eq:eigenvalue-problem-2},
we need the formula for the action of $T_{i}$~\cite{F}:
\begin{gather}
T_{i}E_{\mu}=\left\{
\begin{array}{ll}
a_{i}(\mu)E_{\mu}+b_{i}(\mu)E_{s_{i}\mu} & (\mu_{i}>\mu_{i+1}), \\
E_{\mu} & (\mu_{i}=\mu_{i+1}), \\
a_{i}(\mu)E_{\mu}+E_{s_{i}\mu} & (\mu_{i}<\mu_{i+1}).
\end{array}
\right.
\label{eq:action-of-T}
\end{gather}
Here $a_{i}(\mu)$ and $b_{i}(\mu)$ are def\/ined by
\begin{gather*}
a_{i}(\mu):=\frac{1}{c_{i}(\mu)-c_{i+1}(\mu)}, \qquad
b_{i}(\mu):=1-a_{i}(\mu)^{2}.
\end{gather*}
The formula \eqref{eq:action-of-T} can be obtained from the corresponding formula for
the non-symmetric Jack polynomials and the isomorphism $\Phi$ \eqref{eq:inversion-generic}.

Now let us check \eqref{eq:eigenvalue-problem-2}.
We abbreviate $\lambda(k,s,m,l,r)$ to $\lambda$ in the rest of the proof.
If $\delta_{i}=\delta_{i+1}$, we have $\lambda_{i}>\lambda_{i+1}$,
$a_{i}(\lambda)=-1$ and $b_{i}(\lambda)=0$.
Hence we get $T_{i}E_{\lambda}=-E_{\lambda}$
from \eqref{eq:action-of-T}.
\end{proof}

\subsection[Factorized formula for the case of $(k,r,s)=(N,2,0)$]{Factorized formula for the case of $\boldsymbol{(k,r,s)=(N,2,0)}$}\label{sec:factorized}

In \cite{T} a new integral formula for solutions of the rational qKZ equation
of level zero is constructed.
In the construction we use as an ingredient
a special solution of level one given as follows.
We consider the case of $n=Nm$ for some $m \ge 1$,
and $d_{1}=\cdots= d_{N}=m$.
Then the special solution is given in the form $G=KF$, where
$K=K(z_{1}, \ldots , z_{n})$ is the function def\/ined by~\eqref{eq:K-Gamma},
and
\begin{gather*}
F_{\delta}=\prod_{1\le i<j \le n \atop \delta_{i}=\delta_{j}}(z_{i}-z_{j}-1).
\end{gather*}
It is obtained as a limit of a matrix element
of the vertex operators associated with the quantum af\/f\/ine algebra
$U_{q}(\widehat{sl}_{N})$.
It is also found by Di Francesco and Zinn-Justin in \cite{DZ}.

The special solution above is contained in
the polynomial solutions given in Proposition~\ref{prop:polynomial-solution} as follows.
\begin{proposition}\label{prop:factorized-formula}
Let $(k, r)=(N, 2)$.
For $n=Nm+l$ $(m \ge 1, 0\le l<N)$, consider
\begin{gather*}
\lambda(N,0,m,l,2)=( (m-1, \ldots , 1, 0)^{N-l}, (m, \ldots , 1, 0)^{l}).
\end{gather*}
Then we have
\begin{gather}
\xi_{\lambda(N,0,m,l,2)}(x_{1}, \ldots , x_{n})\big|_{\kappa=-(N+1)}=
\prod_{1\le i<j \le n \atop \delta_{i}=\delta_{j}}(x_{i}-x_{j}),
\label{eq:factorized-xi} \\
E_{\lambda(N,0,m,l,2)}(z_{1}, \ldots , z_{n})\big|_{\kappa=-(N+1)}=
\prod_{1\le i<j \le n \atop \delta_{i}=\delta_{j}}(z_{i}-z_{j}-1),
\label{eq:factorized-E}
\end{gather}
where $\delta=(1^{m}, \ldots , (N-l)^{m}, (N-l+1)^{m+1}, \ldots , N^{m+1})$.
\end{proposition}

\begin{proof}
Hereafter we set $\kappa=-(N+1)$ and abbreviate $\lambda_{0}=\lambda(N,0,m,l,2)$.
Denote the right hand sides of \eqref{eq:factorized-xi} and \eqref{eq:factorized-E}
by $\tilde{\xi}(x)$ and $\tilde{E}(z)$, respectively.

In \cite[Proposition 5.4]{KT} the following formula is proved:
\begin{gather*}
E_{\lambda_{0}}^{q, t}(X)|_{q=t^{-(N+1)}}=
\prod_{1\le i<j \le n \atop \delta_{i}=\delta_{j}}(X_{i}-t^{-1}X_{j}),
\end{gather*}
where $E_{\lambda_{0}}^{q,t}$ is the non-symmetric Macdonald polynomial.
Applying the substitution \eqref{eq:limit}, we see that the right hand side above
has the expansion
\begin{gather*}
\epsilon^{(N-l)\binom{m}{2}+l\binom{m+1}{2}}\big( \tilde{E}(z)+o(1) \big)
\end{gather*}
at $\epsilon=0$.
{}From the argument in the end of Section~\ref{sec:shifted-singular},
$\tilde{E}$ is shifted singular, and hence
$\tilde{\xi}=\Phi^{-1}(\tilde{E})$ is singular.
Hence the Specht module $\tilde{F}:=\mathbb{C}\mathfrak{S}_{n}\cdot\tilde{\xi}$
of isotype $(N^{m}, l)$ consists of singular polynomials.
{}From the result of \cite{D1, D2} (see Section 6 of \cite{D2}),
$\widetilde{F}$ has the basis $\{ \xi_{\mu} \, | \, \mu \in \Lambda^{\rm rev}(\lambda_{0}^{+})\}$,
where $\Lambda^{\rm rev}(\lambda_{0}^{+})$ is the set of reverse lattice permutations of $\lambda_{0}^{+}$.
Since $\lambda_{0}$ is minimal in $\Lambda^{\rm rev}(\lambda_{0}^{+})$
with respect to the dominance order
and $\tilde{\xi}(x)=x^{\lambda_{0}}+(\hbox{lower order terms})$,
we obtain $\tilde{\xi}=\xi_{\lambda_{0}}$.
Sending the both hand sides by $\Phi$, we get $\tilde{E}=E_{\lambda_{0}}$.
\end{proof}

\subsection*{Acknowledgements}

Research of SK is supported by Grant-in-Aid for Scientif\/ic Research (C) No.\,19540228.
Research of MN is supported by Grant-in-Aid for Young Scientists (B) No.\,18749001.
Research of YS is supported by Grant-in-Aid for Scientif\/ic Research (C) No.\,20540009.
Research of YT is supported by Grant-in-Aid for Young Scientists (B) No.\,20740088.

The authors thank Takeshi Suzuki for helpful discussions.
They are also deeply grateful to Charles Dunkl for valuable comments
and pointing out an error in Section~\ref{sec:special-solutions} of the earlier draft.
YT thanks Takeshi Ikeda, Masahiro Kasatani and Hiroshi Naruse for discussions.

\pdfbookmark[1]{References}{ref}
\LastPageEnding


\begin{thebibliography}{99}

\footnotesize\itemsep=-0.5pt

\bibitem{C1}
Cherednik I.,
Nonsymmetric Macdonald polynomials,
{\it Int. Math. Res. Not.} {\bf 1995} (1995), no.~10, 483--515.

\bibitem{C}
Cherednik I.,
Double af\/f\/ine Hecke algebras, Knizhnik--Zamolodchikov equations, and Macdonald's operators,
{\it Int. Math. Res. Not.} {\bf 1992} (1992), no.~9, 171--180.

\bibitem{D1}
Dunkl C.F.,
Singular polynomials for the symmetric groups,
{\it Int. Math. Res. Not.} {\bf 2004} (2004),  no.~67, 3607--3635, \href{http://arxiv.org/abs/math.RT/0403277}{math.RT/0403277}.

\bibitem{D2}
Dunkl C.F.,
Singular polynomials and modules for the symmetric groups,
{\it Int. Math. Res. Not.} {\bf 2005} (2005),  no.~39, 2409--2436, \href{http://arxiv.org/abs/math.RT/0501494}{math.RT/0501494}.

\bibitem{DDO}
Dunkl C.F., de Jeu M.F.E., Opdam E.M.,
Singular polynomials for f\/inite ref\/lection groups,
{\it Trans. Amer. Math. Soc.} {\bf 346} (1994),  237--256.

\bibitem{DZ}
Di Francesco P., Zinn-Justin P.,
Quantum Knizhnik--Zamolodchikov equation, generalized Razumov--Stroganov sum rules and
extended Joseph polynomials,
{\it J. Phys. A: Math. Gen.} {\bf 38} (2005),   L815--L822, \href{http://arxiv.org/abs/math-ph/0508059}{math-ph/0508059}.

\bibitem{F}
Forrester P., Log-gases and Random matrices, in progress.

\bibitem{FJMM}
Feigin B., Jimbo M., Miwa T., Mukhin E.,
A dif\/ferential ideal of symmetric polynomials spanned by Jack polynomials at $\beta=-(r-1)/(k+1)$,
{\it Int. Math. Res. Not.} {\bf 2002} (2002), no.~23, 1223--1237, \href{http://arxiv.org/abs/math.QA/0112127}{math.QA/0112127}.

\bibitem{FR}
Frenkel I.B., Reshetikhin N.Yu.,
Quantum af\/f\/ine algebras and holonomic dif\/ference equations,
{\it Comm. Math. Phys.} {\bf 146} (1992),   1--60.

\bibitem{J}
Joseph A.,
On the variety of a highest weight module,
{\it J. Algebra} {\bf 88} (1984), 238--278.

\bibitem{Knop}
Knop F.,
Symmetric and non-symmetric quantum Capelli polynomials,
{\it Comment. Math. Helv.} {\bf 72} (1997), 84--100, \href{http://arxiv.org/abs/q-alg/9603028}{q-alg/9603028}.

\bibitem{KT}
Kasatani M., Takeyama Y.,
The quantum Knizhnik--Zamolodchikov equation and non-symmetric Macdonald polynomials,
{\it Funkcial. Ekvac.} {\bf 50} (2007),  491--509, \href{http://arxiv.org/abs/math.QA/0608773}{math.QA/0608773}.

\bibitem{M}
Macdonald I.G.,
Af\/f\/ine Hecke algebras and orthogonal polynomials,
{\it Ast\'erisque} {\bf 237} (1996), 189--207.

\bibitem{O}
Opdam E.M.,
Harmonic analysis for certain representations of graded Hecke algebras,
{\it Acta Math.} {\bf 175} (1995), 75--121.

\bibitem{Suzuki}
Suzuki T.,
Rational and trigonometric degeneration of the double af\/f\/ine Hecke algebra of type $A$,
{\it Int. Math. Res. Not.} {\bf 2005} (2005),  no.~37, 2249--2262, \href{http://arxiv.org/abs/math.RT/0502534}{math.RT/0502534}.

\bibitem{T}
Takeyama Y.,
Form factors of $SU(N)$ invariant Thirring model,
{\it Publ. Res. Inst. Math. Sci.} {\bf 39} (2003),  59--116, \href{http://arxiv.org/abs/math-ph/0112025}{math-ph/0112025}.

\end{thebibliography}
\end{document}